 \documentclass[11pt,a4paper,reqno]{amsart} 
\usepackage{amsmath}
\usepackage{amsfonts}

\def\be{\begin{equation}}
\def\ee{\end{equation}}
\def\bea{\begin{eqnarray}}
\def\eea{\end{eqnarray}}
\def\bea*{\begin{eqnarray*}}
\def\eea*{\end{eqnarray*}}

\def\by{{\bf{y}}}
\def\bp{{\bf{p}}}
\def\bq{{\bf{q}}}
\def\br{{\bf{r}}}
\def\bs{{\bf{s}}}
\def\ba{\bf{a}}
\def\bb{\bf{b}}

\begin{document}\date{13/6/19}
\title{An example of the geometry of a 5th-order ODE: the metric on the space of conics in ${\mathbb{CP}}^2$}
\author{Maciej Dunajski}
\address{Department of Applied Mathematics and Theoretical Physics\\ 
University of Cambridge\\ Wilberforce Road, Cambridge CB3 0WA, UK.}
\email{m.dunajski@damtp.cam.ac.uk}

\author{Paul Tod}
\address{The Mathematical Institute\\
Oxford University\\
Woodstock Road, Oxford OX2 6GG\\ UK.
}
\email{tod@maths.ox.ac.uk}

\maketitle
\begin{abstract}
As an application of the method of \cite{dt05}, we find the metric and connection on the space of conics in $\mathbb{CP}^2$ determined as the solution space of the ODE (\ref{1}). These calculations underpin the twistor construction of the Radon transform on conics in $\mathbb{CP}^2$ described in \cite{dt16}. Two further examples of the method are provided.

\end{abstract}
\section{General theory}
In \cite{dt05}, a particular fifth-order ODE whose solutions are the conics in $\mathbb{CP}^2$ was noted as an example for which the W\"unschmann conditions were satisfied, and a torsion-free, $GL(2)$ (or paraconformal) 
connection exists on the moduli 
space $M$ of solutions, while the ODE is not contact equivalent to the trivial fifth-order equation $y^{(5)}=0$. (The fact that this ODE has solutions which are these conics goes back at least to Halphen \cite{h} and may go back to 
Monge \cite{sylvester}.) In this note we spell out all the steps to finding the metric and curvature properties of $M$. These 
calculations are interesting in their own right, for illustrating the method, and they crucially underpin  the twistor construction of the Radon transform on conics in $\mathbb{CP}^2$ described in \cite{dt16}.

\medskip

We begin then with the fifth-order ODE
\be\label{1}
y^{(5)}=\Lambda(x,y,p,q,r,s)=-\frac{40}{9}\frac{r^3}{q^2}+5\frac{rs}{q}
\ee
where $p=y'$ (which doesn't appear yet), $q=y'',r=y'''$ and $s=y''''$. With the conventions of \cite{dt05} we calculate the partial derivatives of $\Lambda$ as
\[\Lambda_x=0=\Lambda_0=\Lambda_1, \Lambda_2=\frac{80}{9}\frac{r^3}{q^3}-5\frac{rs}{q^2},\Lambda_3=5\frac{s}{q}-\frac{40}{3}\frac{r^2}{q^2},\Lambda_4=5\frac{r}{q},\]
and it is straightforward to verify that the W\"unschmann conditions, as in \cite{dt05}, are satisfied so that the moduli space of solutions $M$ admits a torsion-free, $GL(2)$ (or paraconformal) 
connection, defined from $\Lambda$ and its derivatives.

Write the solution as
\[y=Z(x,X^{\bf a}),{\bf a}=1\ldots 5,\]
where $X^{\bf a}$ are coordinates on $M$, and concrete indices are bold. It will eventually be convenient to use $(y,p,q,r,s)$ at some fixed $x$ for $X^{\bf{a}}$, when we'll write them 
$(\by,\bp,\bq,\br,\bs)$.

To say that $M$ has a torsion-free, $GL(2)$-connection \cite{dt05} is to say that the tangent bundle is a symmetric fourth paper of a complex rank-2 spinor bundle, with compatible torsion-free connection preserving 
the spinor symplectic form. Following the method of \cite{dt05}, we impose
\[\nabla_aZ:=Z_{,a}=\iota_A\iota_B\iota_C\iota_D,\]
where $\iota_A$ is a spinor field (and abstract indices are italic). With a slight change from \cite{dt05} we shall suppose
\[\iota'_A=Po_A\]
where prime means $d/dx$ and $P$ is to be found, and $o_A\iota^A=1$ or equivalently the spinor symplectic form is
\[\epsilon_{AB}:=o_A\iota_B-o_B\iota_A\]
and is independent of $x$.

This entails
\[o'_A=Q\iota_A,\]
with $Q$ also to be  found. Compress notation by introducing a constant spinor $\alpha^A$ and writing $\iota=\iota_A\alpha^A,o=o_A\alpha^A$. Then with $t^a=\alpha^A\alpha^B\alpha^C\alpha^D$ write $d\by=t^aZ_a$ etc. Then
\[d\by=\iota^4\]
\[d\bp=(d\by)'=(Z_at^a)'=(\iota^4)'=4Po\iota^3,\]
\[d\bq=(d\bp)'=4PQ\iota^4+4P'o\iota^3+12P^2o^2\iota^2,\]
and
\[d\br=(d\bq)'=A\iota^4+Bo\iota^3+Co^2\iota^2+24P^3o^3\iota\]
with $A,B,C,D$ to be found. By differentiating we obtain
\[A=8P'Q+4PQ',B=4P''+40P^2Q,C=36PP'.\]
Next
\[d\bs=E\iota^4+Fo\iota^3+Go^2\iota^2+Ho^3\iota+24P^4o^4\]
with $E,F,G,H$ to be found, and again by differentiating one calculates
\[E=A'+BQ,F=B'+4PA+2QC,G=C'+3PB+72P^3Q,H=144P^2P'.\]
At the next stage, from (\ref{1}) we have
\[Z^{(5)}_a=\Lambda_2Z''_a+\Lambda_3Z'''_a+\Lambda_4Z''''_a\]
so that
\[(d\bs)'=t^aZ^{(5)}_a=\Lambda_2d\bq+\Lambda_3d\br+\Lambda_4d\bs.\]
Calculating the left-hand-side and equating coefficients gives:
\begin{itemize}
\item
From $o^4$ a differential equation for $P$:
\[96P^3P'+HP=24P^4\Lambda_4,\]
whence
\[240P^3P'=120P^4\frac{\br}{\bq}=120P^4\frac{Z'''}{Z''}\]
which integrates to give
\[P=(Z'')^{1/2}=\bq^{1/2}.\]
\item
From $o^3\iota$ an algebraic equation for $Q$ which solves as

\[Q=\frac{1}{48}\frac{(Z''')^2}{(Z'')^{5/2}}=\frac{1}{48}\frac{\br^2}{\bq^{5/2}}.\]

\item There should then be three identities from the remaining three terms, but we defer considering them for a moment.
\end{itemize}
This choice of $P$ and $Q$ imply
\[A=-\frac{\br^3}{8\bq^3}+\frac{\br\bs}{6\bq^2},\;B=\frac{2\bs}{\bq^{1/2}}-\frac{\br^2}{6\bq^{3/2}},\;C=18\br,\]
and then
\[E=\frac{\bs^2}{6\bq^2}+\frac{\br^2\bs}{6\bq^3}-\frac{319}{864}\frac{\br^4}{\bq^4},\;F=\frac{28}{3}\frac{\br\bs}{\bq^{3/2}}-\frac{151}{18}\frac{\br^3}{\bq^{5/2}},\]
\[G=24\bs+\frac{\br^2}{\bq},\;H=72\bq^{1/2}\br.\]
Now it is straightforward to check that the three identities hold.

\medskip

Note the inverse relation between the basis defined by the spinor dyad (call this the spinor pentad) and the coordinate basis:
\[e^1:=\iota^4=d\by\]
\[e^2:=o\iota^3=\frac{1}{4P}d\bp\]
\[e^3:=o^2\iota^2=\frac{1}{12\bq}(d\bq-\frac{\br}{2\bq}d\bp-\frac{\br^2}{12\bq^2}d\by)\]
\[e^4:=o^3\iota=\frac{1}{24\bq P}\left(d\br-\frac{3\br}{2\bq}d\bq+(\frac{19\br^2}{24\bq^2}-\frac{\bs}{2\bq})d\bp+(\frac{\br^3}{4\bq^3}-\frac{\br\bs}{6\bq^2})d\by\right)\]
\[e^5:=o^4=\frac{1}{24\bq^2}\left(d\bs-\frac{3\br}{\bq}d\br+(-\frac{2\bs}{\bq}+\frac{53\br^2}{12\bq^2})d\bq+(\frac{\br\bs}{6\bq^2}-\frac{17\br^3}{72\bq^3})d\bp+(-\frac{\bs^2}{6\bq^2}+\frac{\br^2\bs}{2\bq^3}-\frac{323\br^4}{864\bq^4})d\by\right)\]
By a general argument the metric may be obtained from the spinor symplectic form as
\[g_{ABCD.PQRS}=\delta_{(A}^{\;K}\delta_B^{\;L}\delta_C^{\;M}\delta_{D)}^{\;N}\epsilon_{KP}\epsilon_{LQ}\epsilon_{MR}\epsilon_{NS},\]
which is equivalent to
\be\label{met66} g=2e^1\odot e^5-8e^2\odot e^4+6e^3\odot e^3,\ee
but in the next section we follow a different route.

\section{Metric}
We obtain the contravariant metric in the chosen coordinates by starting from the condition $g^{ab}Z_aZ_b=0$ and differentiating repeatedly. We use the expressions for the coordinate basis in terms of the spinor pentad to find at once
\[0=g(d\by,d\by)=g(d\by,d\bp)=g(d\by,d\bq)=g(d\by,d\br)=g(d\bp,d\bp)=g(d\bp,d\bq),\]
then
\[g(d\by,d\bs)=Z^aZ''''_a=24P^4\mbox{  so  }g(d\bp,d\br)=Z'^aZ'''_a=-24P^4\]
and then 
\[g(d\bq,d\bq)=Z''^aZ''_a=24P^4,\]
where recall that $P^2=\bq$.

Next
\[g(d\bp,d\bs)=Z'^aZ''''_a=-144P^3P'\mbox{  so  }g(d\bq,d\br)=Z''^aZ'''_a=48P^3P',\]
and recall $2P^3P'=\bq\br$.

For the rest a harder calculation gives
\[g(d\bq,d\bs)=48\bq\bs-32\br^2,\;g(d\br,d\br)=56\br^2-24\bq\bs,\] and finally \[g(d\br,d\bs)=\frac{160}{3}\frac{\br^3}{\bq}-16\br\bs,\;g(d\bs,d\bs)=104\bs^2-320\frac{\br^2\bs}{\bq}+\frac{2560}{9}\frac{\br^4}{\bq^2}.\]
The contravariant metric is then

  \[
(g^{{\bf ab}})=\left(\begin{array}{ccccc}
           0 & 0 & 0 & 0 & 24\bq^2\\
          0 & 0 & 0 & -24\bq^2 & -72\bq\br\\
           0 & 0 & 24\bq^2 & 24\bq\br & 48\bq\bs-32\br^2 \\
           0 &-24\bq^2 & 24\bq\br &56\br^2-24\bq\bs & \frac{160}{3}\frac{\br^3}{\bq}-16\br\bs\\
           24\bq^2 &-72\bq\br  &48\bq\bs-32\br^2&\frac{160}{3}\frac{\br^3}{\bq}-16\br\bs&104\bs^2-320\frac{\br^2\bs}{\bq}+\frac{2560}{9}\frac{\br^4}{\bq^2}\\
\end{array}\right),
\]
with coordinates ordered $(\by,\bp,\bq,\br,\bs)$.

This inverts to:

  \[
(g_{{\bf ab}})=\left(\begin{array}{ccccc}
           \frac{\br^2\bs}{24\bq^5}-\frac{5\br^4}{162\bq^6}-\frac{\bs^2}{72\bq^4} & \frac{\br\bs}{72\bq^4}-\frac{\br^3}{54\bq^5} & \frac{13}{72}\frac{\br^2}{\bq^4}-\frac{s}{12\bq^3} & -\frac{\br}{8\bq^3} & (24\bq^2)^{-1}\\
          * & \frac{\bs}{24\bq^3}-\frac{\br^2}{18\bq^4} & \frac{\br}{24\bq^3} & -(24\bq^2)^{-1} & 0\\
           * & * & (24\bq^2)^{-1} & 0 & 0 \\
           * & * & 0 &0 & 0\\
           * &0  &0&0&0\\
\end{array}\right),
\]
where $*$ indicates a term fixed by symmetry.
\begin{itemize}
\item
Now we can check directly that
\[(g_{\bf ab}dX^{\bf a}dX^{\bf b})'=0,\]
so that the metric is independent of the fixed but arbitrary choice we made of $x$, and also check that it agrees with (\ref{met66}).
\item Note also that $g_{\bf yy}=\frac{\br^2\bs}{24\bq^5}-\frac{5\br^4}{162\bq^6}-\frac{\bs^2}{72\bq^4}$ is constant: this is a first-integral of (\ref{1}).
\item A Maple calculation shows that the metric is Einstein with scalar curvature $R=-60$. That this is Einstein is to be expected from the general theory of symmetric spaces, as the space of conics is the symmetric space $SL(3, \mathbb{C})/SO(3,\mathbb{C})$. 
This space contains two real forms: the positive definite $SL(3,\mathbb{R})/SO(3, \mathbb{R})$, and pseudo--Riemanian $SL(3, \mathbb{R})/SO(1, 2)$, where the
metric has signature $(2, 3)$. This later case is what we obtain if (\ref{1}) is regarded as a real ODE, and $(x, y)$ are taken to be real. In \cite{dt16} the
Riemannian form was used.
\item
The Laplacian has no first-order derivatives and so is just
\[\Delta=g^{ab}\partial_a\partial_b.\]
Equivalently, the coordinates are harmonic.
\end{itemize}

\section{Connection}
We assume the Levi-Civita derivative $\nabla_a$ extends to spinors so the derivative of the spinor dyad can be written
\[\nabla_ao_B=\phi_ao_B+\psi_a\iota_B,\;\;\nabla_a\iota_B=\chi_ao_B+\lambda_a\iota_B,\]
for vectors $\phi_a,\psi_a,\chi_a,\lambda_a$ to be found, and that this preserves $\epsilon_{AB}$:
\[0=\nabla_a\epsilon_{AB}=\nabla_a(o_A\iota_B-o_B\iota_A).\]
Therefore $\phi_a+\lambda_a=0$ so we may eliminate $\lambda_a$.
Note that
\[\phi_a=\iota^B\nabla_ao_B=o^B\nabla_a\iota_B,\;\;\psi_a=-o^B\nabla_ao_B,\;\;\chi_a=\iota^B\nabla_a\iota_B.\]

We also want this extended derivative to be torsion-free. Recall
\[\by_a=\iota_A\iota_B\iota_C\iota_D,\;\;\bp_a=4\bq^{1/2}o_{(A}\iota_B\iota_C\iota_{D)},\]
so that
\[\nabla_a\by_b=\bq^{-1/2}\chi_a\bp_b-4\phi_a\by_b.\]
Torsion-free-ness necessarily requires $\nabla_{[a}y_{b]}=0$ and this will also be sufficient as priming it up to four times shows. Thus
\[\chi_{[a}p_{b]}-4\bq^{1/2}\phi_{[a}y_{b]}=0\] from which it follows that
\[
\phi_a=\alpha \by_a-\bq^{-1/2}\gamma \bp_a \]
\[\chi_a=4\gamma \by_a+\delta \bp_a\]
for some $\alpha,\gamma,\delta$ to be found.

We may calculate primes of $\chi_a,\phi_a,\psi_a$ by using $\iota'_A=Po_A,o'_A=Q\iota_A$ to obtain:
\[\chi'_a=(\iota^B\nabla_a\iota_B)'=2P\phi_a+\nabla_aP\]
\[\psi'_a=-2Q\phi_a+\nabla_aQ\]
\[\phi'_a=-P\psi_a+Q\chi_a,\]
where, recall, $P=\bq^{1/2},Q=\frac{1}{48}\frac{\br^2}{\bq^{5/2}}$.

Substitute into $\chi'_a$:
\[4\gamma'\by_a+4\gamma \bp_a+\delta'\bp_a+\delta \bq_a=2\bq^{1/2}(\alpha \by_a-\bq^{-1/2}\gamma \bp_a)+\frac{1}{2}\bq^{-1/2}\bq_a,\]
so that
\[\delta=\frac{1}{2}\bq^{-1/2},\]
\[6\gamma+\delta'=0\mbox{  whence  }\gamma=\frac{1}{24}\bq^{-3/2}\br,\]
and
\[\alpha=2\bq^{-1/2}\gamma'=\frac{\bs}{12\bq^2}-\frac{\br^2}{8\bq^3}.\]
Thus $\phi_a$ and $\chi_a$ are now known. For $\psi_a$ consider $\phi'_a$:
\[P\psi_a=-\phi'_a+Q\chi_a=-(\alpha \by_a-\bq^{-1/2}\gamma \bp_a)'+Q(4\gamma \by_a+\delta \bp_a)\]
\[=(-\alpha'+4Q\gamma)\by_a+(-\alpha+(\bq^{-1/2}\gamma)'+Q\delta)\bp_a+\bq^{-1/2}\gamma \bq_a\]
whence
\[\psi_a=-\frac{1}{864}\frac{\br^3}{\bq^{9/2}}\by_a+\left(\frac{5}{96}\frac{\br^2}{\bq^{7/2}}-\frac{\bs}{24\bq^{5/2}}\right)\bp_a+\frac{\br}{24\bq^{5/2}}\bq_a.
\]
The equation for $\psi'_a$ should now be an identity and indeed it is.

\medskip

Note now that

\[\iota^A\iota^B\iota^C\nabla_{ABEF}\iota_C=0,\]
 This is the condition for integrability of the distribution spanned by $\iota^A\nabla_{ABCD}$ and the integral manifolds of the distribution are the surfaces of constant $y$. Such a surface is defined by all 
 conics through a fixed point of ${\mathbb{CP}}^2$, which will recur in the final section. In fact we have here a stronger result:
 
 \be\label{int1}\iota^A\nabla_{ABCD}\iota_E=\iota^A(\chi_{ABCD}o_E-\phi_{ABCD}\iota_E)=\iota_B \iota_C \iota_D(\frac12 o_E-\gamma\iota_E),
 \ee
 a formula which is needed in \cite{dt16}.

\section{The $SO(3)$-structure}
We first recall some $SO(3)$-theory following \cite{BN} and \cite{F}. The metric $g_{ab}$ is defined from the spinor epsilon as in (\ref{met66}) but here we introduce a new notation for this:
\[g_{ae}=g_{ABCDEFGH}={\mathcal{S}}_{(ABCD)}\left(\epsilon_{AE}\epsilon_{BF}\epsilon_{CG}\epsilon_{DH}\right),\]
where the symbol ${\mathcal{S}}_{(ABCD)}$ is introduced to define symmetrisation of the following expression over the indices $ABCD$ with the usual factor $(4!)^{-1}$. 
We may define an analogous symmetric tensor $G_{aep}$ from six epsilons by
\[G_{aep}=G_{ABCDEFGHPQRS}={\mathcal{S}}_{(ABCD)}{\mathcal{S}}_{(EFGH)}\left(\epsilon_{AE}\epsilon_{BF}\epsilon_{GP}\epsilon_{HQ}\epsilon_{CR}\epsilon_{DS}\right).\]
It is straightforward to check that
\[G_{abc}=G_{(abc)},\;\;g^{ab}G_{abc}=0,\;\;\nabla_aG_{bcd}=0,\]
and the normalisation 
\be\label{G0}6G^e_{\;\;a(b}G_{cd)e}=g_{a(b}g_{cd)}\ee
holds. 

More identities follow: trace (\ref{G0}) to obtain
\[G_{efa}G^{ef}_{\;\;\;\;b}=\frac{7}{12}g_{ab}\mbox{   and   }G_{abc}G^{abc}=\frac{35}{12}.\]
Commute derivatives on $G_{abc}$ to obtain a condition on the curvature tensor:
\be\label{G1}R_{abc}^{\;\;\;\;\;\;(d}G^{ef)c}=0.\ee
Define
\[\chi_{abcd}=6G^e_{\;\;ab}G_{cde},\;\;F_{bcad}=\chi_{a[bc]d},\]
and claim
\be\label{G2}\chi_{abcd}=\chi_{(abcd)}+\frac{2}{3}F_{bcad}+\frac{2}{3}F_{bdac},\ee
with
\[\chi_{(abcd)}=6G^e_{\;\;(ab}G_{cd)e}=g_{a(b}g_{cd)}.\]
Expand (\ref{G1}):
\[R_{abc}^{\;\;\;\;\;\;d}G^{efc}+R_{abc}^{\;\;\;\;\;\;e}G^{fdc}+R_{abc}^{\;\;\;\;\;\;f}G^{dec}=0\]
and contract with $G_{efp}$ to deduce
\be\label{G4}R_{abcd}F^{cd}_{\;\;\;\;\;\;pq}=\frac{7}{4}R_{abpq},\ee
after relabelling of indices. 

We need these identities in the next section.
\section{A system of equations}
In \cite{dt16} and following Moraru \cite{m1} we consider the system of equations:
\begin{eqnarray}
G_a^{\;\;bc}\nabla_b\nabla_cF&=&\lambda\nabla_aF,\label{G55}\\
\Delta F:=g^{ab}\nabla_a\nabla_bF&=&\mu F,\label{G6}
\end{eqnarray}
on a scalar $F$, where $\lambda,\mu$ are real constants. These can be written down in any $SO(3)$-structure but we are interested principally in the case of Section~1, which is also Einstein.

Compress notation by writing $F_a=\nabla_aF$ then from (\ref{G55})
\[6\lambda^2F^a=6G^{abc}G_b^{\;\;de}\nabla_c\nabla_dF_e=\chi^{acde}\nabla_c\nabla_dF_e\]
\[=(g^{a(c}g^{de)}+\frac{2}{3}F^{cdae}+\frac{2}{3}F^{cead})\nabla_c\nabla_dF_e\]
using identites from the previous section. Here the first term is
\[\frac{1}{3}(g^{ac}g^{de}+g^{ad}g^{ec}+g^{ae}g^{cd})\nabla_c\nabla_dF_e\]
\[=\frac{1}{3}(\nabla^a\Delta F+2\nabla_c\nabla^aF^c)\]
\[=\frac{1}{3}(\nabla^a\Delta F+2R^{ab}F_b+2\nabla^a\Delta F)\]
\[=\nabla^a(\Delta F+\frac{2}{15}RF),\]
using the Einstein condition.

The other two terms become
\[\frac{4}{3}F^{cdae}\nabla_c\nabla_dF_e=-\frac{2}{3}F^{cdae}R_{cdfe}F^f\]
\[=-\frac{2}{3}.\frac{7}{4}R^{ae}_{\;\;\;\;\;fe}F^f=-\frac{7}{6}R^a_{\;\;\;f}F^f=-\frac{7}{30}RF^a,\]
using the Einstein condition again.

Putting these together
\[6\lambda^2F_a=\nabla_a(\Delta F-\frac{1}{10}RF),\]
whence
\be\label{G7}\mu=6\lambda^2+\frac{R}{10}.\ee
Conversely, a solution $F$ of (\ref{G55}) with some $\lambda$ will necessarily satisfy (\ref{G6}) with the value of $\mu$ given by (\ref{G7}), possibly after adding a constant to $F$.

\medskip

In spinor notation the system (\ref{G55}) can be written
\[\Box_{ABCD}F:=\nabla_{(AB}^{\;\;\;\;\;\;\;\;EF}\nabla_{CD)EF}F=\lambda\nabla_{ABCD}F,\]
accompanied by
\[\Delta F =\mu F.\]

\medskip

To write out the system in coordinates we need to calculate two sets of quantities:
\[G_a^{\;\;bc}\nabla_b\nabla_cX^{\bf{a}}\mbox{   and   }G_a^{\;\;bc}\nabla_bX^{\bf{b}}\nabla_cX^{\bf{c}},\]
but we have all the necessary information for these, so we may assume them known.

For a function $F$, the one-form $G_a^{\;\;bc}\nabla_b\nabla_cF$ decomposes in the coordinate basis as:
\[d\by: 4\bq F_{\by\bq}+6\br F_{\by\br}+8\bs F_{\by\bs}-2\bq F_{\bp\bp}-2\br F_{\bp\bq}-2\bs F_{\bp\br}-2\bs' F_{\bp\bs}-F_\by\]
\[d\bp: 6\bq F_{\by\br}+16\br F_{\by\bs}-2\bq F_{\bp\bq}-4\br F_{\bp\br}-6\bs F_{\bp\bs}-F_\bp\]
\[d\bq: 4\bq F_{\by\bs}+2\bq F_{\bp\br}-2\br F_{\bq\br}-2\bq F_{\bq\bq}+(-16\bs+\frac{80\br^2}{3\bq})F_{\bq\bs}
+(7\bs-\frac{40\br^3}{3\bq})F_{\br\br}\]\[+(\frac{70\br\bs}{3\bq}-\frac{400\br^3}{9\bq^2})F_{\br\bs}+ (-\frac{70}{3}\frac{\bs^2}{\bq}+\frac{320}{3}\frac{\br^2\bs}{\bq^2}-\frac{3200}{27}\frac{\br^4}{\bq^3})F_{\bs\bs}-F_\bq\]
\[d\br: 4\bq F_{\bp\bs}-16\br F_{\bq\bs}-2\bq F_{\bq\br} +6\br F_{\br\br}+(-2\bs+\frac{80\br^2}{3\bq})F_{\br\bs}+
(-\frac{80}{3}\frac{\br\bs}{\bq}+\frac{640}{9}\frac{\br^3}{\bq^2}) F_{\bs\bs}-F_\br\]
\[d\bs: 4\bq F_{\bq\bs}-3\bq F_{\br\br}-12\br F_{\br\bs}+(8\bs-\frac{80\br^2}{3\bq})F_{\bs\bs}-F_\bs,\]
and this must be equated to $\lambda dF$.

This system is considered further in \cite{dt16} and it is shown there that solutions are given as follows: pick $f(x,y)$ and perform the integral
\[ F(X^{\bf a})=\int f(x,Z(x,X^{\bf a}))q^{1/3}dx\]
over a suitable contour. This is a translation of a formula in \cite{m1} and generates solutions of the system (\ref{G55}-\ref{G6}).

\section{Further examples}
The methods of this paper can be extended to a wider selection of examples but it follows from \cite{dt05} and \cite{gn} that, while there are other fifth-order ODEs giving rise to $SO(3)$-structures in the sense used here, the connection preserving the tensor $G_{abc}$ in general has torsion -- the unique non-trivial torsion-free case is the one presented above. We'll give below an example of another fifth-order ODE leading to an $SO(3)$-structure, and also an example of a fourth-order ODE where the moduli space admits one of Bryant's exotic ${\mathcal{G}}_3$-holonomy connections (\cite{br1}; for this example, the theory in the form we need it can be found in \cite{dt05}).

\subsection{The Fifth-order ODE}

From \cite{gn} and with the notation of (\ref{1}) we consider the equation
\be\label{f1}
y^{(5)}=\Lambda(x,y,p,q,r,s)=\frac{5}{3}\frac{s^2}{r}.\ee

This is readily solved to give 
\[y=c_5+c_4x+c_3x^2+(c_1+c_2x)^{3/2},\]
but the interest in the equation for us is that the relevant W\"unschmann invariants
vanish \cite{gn}. As before we write the solution as $y=Z(x;X^{\ba})$ with conventions for the coordinates $(\by,\bp,\bq,\br,\bs)$ on the moduli space as before, and we introduce spinors with
 \[y_a=\iota_A\iota_B\iota_C\iota_D \mbox{ or  }d\by=(\iota)^4.\]
 We assume
 \[\iota'_A=Po_A,\;\;\;o'_A=Q\iota_A,\]
where $o_A$ forms a normalised spinor dyad with $\iota_A$, and $P,Q$ are to be found.

Then
\[d\bp=(d\by)'=4Po\iota^3,\]
\[d\bq=(d\bp)'=4PQ\iota^4+4P'o\iota^3+12P^2o^2\iota^2\]
\[d\br=(d\bq)'=A\iota^4+Bo\iota^3+Co^2\iota^2+24P^3o^3\iota\]
with
\[A=8P'Q+4PQ',B=4P''+40P^2Q,C=36PP',\]
and so 
\[d\bs=E\iota^4+Fo\iota^3+Go^2\iota^2+Ho^3\iota+24P^4o^4\]
where one calculates
\[E=A'+BQ,F=B'+4PA+2QC,G=C'+3PB+72P^3Q,H=144P^2P'.\]
Finally, using (\ref{f1}),
\[(d\bs)'=\Lambda_\br d\br+\Lambda_\bs d\bs=-\frac{5\bs^2}{3\br^2}d\br+\frac{10\bs}{3\br}d\bs.\]
From the coefficient of $d\bs$:
\[240P^3P'=\frac{10\bs}{3\br}.24P^4\mbox{  whence  }\frac{P'}{P}=\frac{\bs}{3\br}\mbox{  and  }P=\br^{1/3}.\]
Next from the coefficient of $d\br$ we find that $Q=0$, and then the three equations from $d\by,d\bp,d\bq$ are all identities. Summarising:
\[A=0,\;B=\frac{4\bs^2}{3\br^{5/3}},\;C=\frac{12\bs}{\br^{1/3}},\;E=0,\;F=\frac{20\bs^3}{9\br^{8/3}},\;G=\frac{20\bs^2}{\br^{4/3}},\;H=48\bs.\]
The orthonormal basis is
\[e_1=\iota^4=d\by,\]
\[e_2=o\iota^3=\frac{1}{4P}d\bp,\]
\[e_3=o^2\iota^2=\frac{1}{12P^2}(d\bq-\frac{\bs}{3\br}d\bp),\]
\[e_4=o^3\iota=\frac{1}{24P^3}(d\br-\frac{\bs}{\br}d\bq).\]
\[e_5=o^4=\frac{1}{24P^4}(d\bs-\frac{2\bs}{\br}d\br+\frac{\bs^2}{3\br^2}d\bq),\]
with duals
\[E_1=\partial_\by\]
\[E_2=4P(\partial_\bp+\frac{\bs}{3\br}\partial_\bq)\]
\[E_3=12P^2(\partial_\bq+\frac{\bs}{\br}\partial_\br+\frac{5\bs^2}{3\br^2}\partial_\bs)\]
\[E_4=24P^3(\partial_\br+\frac{2\bs}{\br}\partial_\bs)\]
\[E_5=24P^4\partial_s.\]
Now the metric from (\ref{met66}) is
\[g=2(e_1\odot e_5-4e_2\odot e_4+3e_3\odot e_3)\]
\[=\frac{1}{24\br^{4/3}}(2d\by(d\bs-\frac{2\bs}{\br}d\br+\frac{\bs^2}{3\br^2}d\bq  )-2d\bp(d\br-\frac{\bs}{\br}d\bq )+(d\bq-\frac{\bs}{3\br}d\bp)^2),\]
or as a matrix
  \[(g_{\ba\bb})=\left(\begin{array}{ccccc}
0 & 0 & \frac{\bs^2}{72\br^{10/3}} & -\frac{\bs}{12\br^{7/3}} & \frac{1}{24\br^{4/3}}\\
0 & \frac{\bs^2}{216\br^{10/3}}&\frac{\bs}{36\br^{7/3}} & -\frac{1}{24\br^{4/3}} & 0\\
\frac{\bs^2}{72\br^{10/3}} & \frac{\bs}{36\br^{7/3}} & \frac{1}{24\br^{4/3}} & 0 & 0 \\
-\frac{\bs}{12\br^{7/3}} & -\frac{1}{24\br^{4/3}} & 0 &0 & 0\\
\frac{1}{24\br^{4/3}} & 0 &0& 0 &0\\
\end{array}\right).\]
The method of Section 2 to obtain the metric starts from
\[g(d\by,d\by)=0\]
whence by differentiating
\[0=g(d\by,d\bp)=g(d\by,d\bq)=g(d\by,d\br),\;\;g(d\by,d\bs)=24\br^{4/3},\]
\[g(d\bp,d\bp)=0=d(d\bp,d\bq),\;\;g(d\bp,d\br)=-24\br^{4/3},\;\;g(d\bp,d\bs)=-48\br^{1/3}\bs,\]
and so on, culminating in
\[(g^{\ba\bb})=\left(\begin{array}{ccccc}
0 & 0 & 0 & 0 & 24\br^{4/3}\\
0 & 0 & 0 & -24\br^{4/3} & -48\br^{1/3}\bs\\
0 & 0 & 24\br^{4/3} & 16\br^{1/3}s & 24\br^{-2/3}\bs^2 \\
0 & -24\br^{4/3} & 16\br^{1/3}\bs & 8\br^{-2/3}\bs^2 & -\frac{32}{3}\br^{-5/3}\bs^3\\
24\br^{4/3} &-48\br^{1/3}s  & 24\br^{-2/3}\bs^2 & -\frac{32}{3}\br^{-5/3}\bs^3 & \frac{40}{3}\br^{-8/3}\bs^4\\
\end{array}\right).\]
It is straightforward to check that these matrices are inverses, and that the metric has six independent Killing vectors and a homothety, and is scalar-flat but not Ricci-flat.

If we next follow the method of Section 3 to seek a torsion-free spinor connection inducing the Levi-Civita connection on vectors and annihilating $\epsilon_{AB}$ we reach a contradiction, since we know from \cite{dt05} and \cite{gn} that the connection preserving the $SO(3,{\mathbb{C}})$-structure has torsion. We shall leave this example here.

\subsection{The Fourth-order ODE}
The association of an exotic ${\mathcal{G}}_3$-holonomy connection in four-dimensions to a fourth-order ODE satisfying certain conditions is due to Bryant \cite{br1}. The conditions are the vanishing of certain W\"unschmann invariants of the ODE as was made explicit in Theorem 1.3 of \cite{dt05}. An example is provided by the ODE determining the conics in $\mathbb{CP}^2$ which pass through a given fixed point. It is straightforward to check that the relevant W\"unschmann invariants do vanish, so that the moduli space admits what was called a paraconformal structure in \cite{dt05} and this is nontrivial, in the sense that the ODE is not contact-equivalent to the trivial equation $y^{(4)}=0$ by a criterion from \cite{doub} quoted in Theorem 3.5 of \cite{dt05}. There will be a connection preserving the paraconformal structure but it will necessarily have torsion. We won't compute it but we will describe it below.

We consider then the ODE satisfied by conics through a fixed point in ${\mathbb{CP}}^2$. These can be taken to have equation
\[ax^2+2bxy+y^2+2cx+2ey=0\]
when the fixed point is $(0,0)$ in an affine patch, and the fourth-order equation annihilating $y(x)$ is
\be\label{f2}
y^{(4)}=\Lambda(x,y,p,q,r)=\frac{4r^2}{3q}+\frac{2xqr+6q^2}{xp-y}-\frac{3x^2q^3}{(xp-y)^2}    .\ee
This looks a little simpler in terms of $W:=xp-y$, when
\[\Lambda=\frac{4r^2}{3q}+\frac{2xqr+6q^2}{W}-\frac{3x^2q^3}{W^2}.\]

Note that for this example we expect a preserved symplectic form or equivalently a symmetric quartic as shown in \cite{dt05}, but not a metric. In coordinates $(\by,\bp,\bq,\br)$ introduced in the now standard way we note that the two-form is 
\be\label{o1}\Omega:=(\iota)^3\wedge(o)^3-3o(\iota)^2\wedge (o)^2\iota\ee
and claim that, in coordinates,
\begin{eqnarray}
\Omega&=&\frac{1}{6\bq^{4/3}(x\bp-\by)}\big(d\by\wedge d\br-d\bp\wedge d\bq+(\frac{4\br}{3\bq}+\frac{2x\bq}{x\bp-\by})d\by\wedge d\bq\nonumber\\
&&-\frac{1}{(x\bp-\by)^2}((x\bp-\by)(x\br-3\bq)-3x^2\bq^2)d\by\wedge d\bp\big). \label{form} 
\end{eqnarray}

To see this, we first introduce spinors as before so with solution $y~=~Z(x;X^{\ba})$ to (\ref{f2}), set
\[d\by=(\iota)^3=\iota_A\iota_B\iota_C.\]
It is important not to confuse $\by$ here with $\by$ in Sections 1-5, which solves (\ref{1}) rather that (\ref{f1}), nor to confuse $\iota_A$ here with $\iota_A$ there.

Suppose
\[\iota'_A=Po_A,\;\;\;o'_A=Q\iota_A\]
with $P,Q$ to be found, then
\[d\bp=(d\by)'=3Po\iota^2,\]
\[d\bq=(d\bp)'=3PQ\iota^3+3P'o\iota^2+6P^2o^2\iota,\]
and
\[d\br=(d\bq)'=A\iota^3+Bo\iota^2+Co^2\iota+Do^3,\]
with
\[A=3PQ'+6P'Q,\]
\[B=3P''+21P^2Q,\]
\[C=18PP'\]
\[D=6P^3.\]
It is convenient to invert the relations between the normalised and coordinate tetrads:
\[(\iota)^3=d\by,\]
\[o(\iota)^2=\frac{1}{3P}d\bp\]
\[o^2\iota=\frac{1}{6P^2}\left(d\bq-\frac{P'}{P}d\bp-3PQd\by\right)\]
\[o^3=\frac{1}{D}\left(d\br-\frac{C}{6P^2}d\bq+(\frac{CP'}{6P^3}-\frac{B}{3P})d\bp+(\frac{CQ}{2P}-A)d\by\right).\]
To write the symplectic form (\ref{o1}) in coordinates we need $P,Q,B$ and $C$.
From the two expressions for $d\bs$:
\[d\bs=d\Lambda =\Lambda_\by d\by+\Lambda_\bp d\bp+\Lambda_\bq d\bq+\Lambda_\br d\br\]
\[=E\iota^3+Fo\iota^2+Go^2\iota+Ho^3\]
we obtain
\[E=A'+BQ,\]
\[F=B'+3AP+2CQ\]
\[G=C'+2BP+3DQ,\]
\[H=D'+CP=36P^2P'.\]

Use the relation of the spinor tetrad to the holonomic tetrad to write:
\[\Lambda_\by d\by+\Lambda_\bp d\bp+\Lambda_\bq d\bq+\Lambda_\br d\br=\Lambda_\by\iota^3+\Lambda_\bp3Po\iota^2\]
\[+\Lambda_\bq(3PQ\iota^3+3P'o\iota^2+6P^2o^2\iota)+\Lambda_\br(A\iota^3+Bo\iota^2+Co^2\iota+Do^3).\]
and then read off corresponding terms. From $o^3$
\[D\Lambda_\br=H=36P^2P'\]
so that 
\[\frac{P'}{P}=\frac16\Lambda_\br=\frac16\left(\frac{8\br}{3\bq}+\frac{2x\bq}{x\bp-\by}\right),\]
and integrate to obtain
\[\log P=\frac{4}{9}\log \bq+\frac13\log(x\bp-\by),\]
dropping the constant of integration. Exponentiating
\[P=\bq^{4/9}(x\bp-\by)^{1/3}.\]
From $o^2\iota$
\[G=C'+2BP+3DQ=C\Lambda_\br+6P^2\Lambda_\bq,\]
i.e.
\[18PP''+18(P')^2+2P(3P''+21P^2Q)+18P^3Q=18PP'\Lambda_\br+6P^2\Lambda_\bq,\]
which solves for $Q$:
\[Q=\frac{1}{60P}\left(-24\frac{P''}{P}-18(\frac{P'}{P})^2+6\Lambda_\bq+18\frac{P'}{P}\Lambda_\br\right),\]
\[=\frac{1}{9P}\left(\frac{2\br^2}{9\bq^2}+\frac{x\br}{3W}-\frac{x^2\bq^2}{W^2}\right).\]
Next
\[D=6P^3=6\bq^{4/3}W,\]
\[C=18PP'=18P^2\left(\frac{4\br}{9\bq}+\frac{x\bq}{3W}\right)=\frac{2\bq^{1/3}}{P}(4\br W+3x\bq^2),\]
\[B=3P''+21P^2Q=P\left(\frac{14\br^2}{\bq^2}+\frac{1}{3W}(16x\br+27\bq)-\frac{7x^2\bq^2}{W^2}\right).\]
Now that we have $P,Q,B$ and $C$ we can be explicit about $\Omega$, substituing for the spinors in (\ref{o1}) and we obtain precisely (\ref{form}).  
We readily check that $\Omega$ is $x$-independent, and closed and non-degenerate in the sense

\[\Omega\wedge\Omega=-\frac{1}{18P^6}d\by\wedge d\bp\wedge d\bq\wedge d\br\neq 0.\]

We could go on to calculate $A$ and check the two remaining equations
\[B'+3AP+2CQ=3P\Lambda_p+3P'\Lambda_q+B\Lambda_r,\]
\[A'+BQ=\Lambda_y+3PQ\Lambda_q+A\Lambda_r,\]
but these must be identities, by general theory.

\medskip

We shall leave this example here but note that it has an interpretation in terms of the five-dimensional example treated in Sections 1-5 above. There the moduli space, say $\mathcal{M}^5$, was the set of all conics in $\mathbb{CP}^2$; here it is the set, say $\mathcal{N}^4$, of such conics through a fixed point. Evidently $\mathcal{N}^4$ is a hypersurface in $\mathcal{M}^5$, and in fact a hypersurface of constant $\by$, using $\by$ in the sense of Section 1. This is a null hypersurface so has only a degenerate metric. The normal to it, using $\iota^A$ from Section 1 is 
\[d\by=(\iota)^4=\iota_A\iota_B\iota_C\iota_D,\]
with $\iota_A$ in the sense of Section 1. A tangent vector to $\mathcal{N}^4$ takes the form
\[V^{ABCD}=\iota^{(A}V^{BCD)}\]
and can be represented by $V^{BCD}$ in $T\mathcal{N}$. A covariant derivative $D_{ABC}$ can be defined on $\mathcal{N}^4$ by
\[D_{ABC}:=\iota^D\nabla_{ABCD},\]
using the Levi-Civita derivative $\nabla_{ABCD}$ from Section 3. Evidently this derivative annihilates functions on $\mathcal{M}^5$ which are constant on $\mathcal{N}^4$, and it preserves $\epsilon_{AB}$ and therefore $\Omega_{ab}$, but it will have torsion as we see by commuting on scalars:
\[(D_{ABC}D_{PQR}-D_{PQR}D_{ABC})f=(\iota^D\nabla_{ABCD}\iota^S)\nabla_{PQRS}f-(\iota^S\nabla_{PQRS}\iota^D)\nabla_{ABCD}f\]
\[=\beta^S(\iota_A\iota_B\iota_C\nabla_{PQRS}-\iota_P\iota_Q\iota_R\nabla_{ABCS})f  \]
where 
\[\beta_A=\frac12o_A-\gamma\iota_A,\]
after substituting from (\ref{int1}), and then
\[=\beta^S\iota_{(A}\iota_B\epsilon_{C)(P}D_{QR)S}f=T_{ap}^{\;\;\;m}D_mf,\]
for a torsion tensor $T_{ap}^{\;\;\;\;\;m}$ which can be expressed in terms of $\epsilon_{AB}$, $\delta_A^{\;B}$ and the vector $\phi_{ABC}=\iota_{(A}\iota_B\beta_{C)}$ as
\[T_{ABC.PQR}^{\;\;\;\;\;\;\;\;\;\;\;\;\;\;\;\;\;LMN}=\epsilon_{E(A}\delta_B^{(L}\phi_{C)(P}^{\;\;\;\;\;\;\;M}\delta_Q^{N)}\delta_{R)}^E.\]

\section*{Acknowledgements}
The work of MD was partially supported by STFC consolidated grant  no. ST/P000681/1. Part of this work was done while PT held the Brenda Ryman Visiting Fellowship in the Sciences at Girton College, Cambridge, and 
he gratefully acknowledges the hospitality of the College and of CMS.

\end{document}